\documentclass[12pt]{article}
\usepackage{latexsym}
\usepackage{amsfonts}
\newtheorem{theorem}{Theorem}
\newtheorem{lemma}{Lemma}
\newtheorem{proposition}{Proposition}

\begin{document}
\def\arge{{\rm large}}
\def\proj{{\rm proj}}
\def\Tr{{\rm Tr}}
\def\grpz{{\mathbb Z}}
\def\grpc{{\mathbb C}}
\def\fbar{{\overline f}}
\def\vcw{{\bf w}}
\def\vcy{{\bf y}}
\def\vcv{{\bf v}}
\def\vcz{{\bf 0}}
\def\vcb{{\bf b}}
\def\vcc{{\bf c}}
\def\vca{{\bf a}}
\def\vcbnmt{{\bf B_{n+m-2}}}
\def\vcbnmo{{\bf B_{n+m-1}}}
\def\vcba{{\bf b_1}}
\def\vcbd{{\bf b_d}}
\def\vcbj{{\bf b_j}}
\def\vcur{{\bf u_r}}
\def\vcer{{\bf e_r}}
\def\vcea{{\bf e_1}}
\def\vced{{\bf e_d}}
\def\vcek{{\bf e_k}}
\def\vcxm{{\bf X_m}}
\def\vcxcm{{\bf X_{cm}}}
\def\vcxtm{{\bf X_{2m}}}
\def\vcxnm{{\bf X_{n+m}}}
\def\vcxno{{\bf X_{n+1}}}
\def\vcxn{{\bf X_n}}
\def\vcbn{{\bf B_n}}
\def\vcxz{{\bf X_0}}
\def\vcbz{{\bf B_0}}
\def\vcbo{{\bf B_1}}
\def\qed{{\hfill $\Box$}}
\baselineskip=20pt
\title{Generating Random Vectors in $(\grpz/p\grpz)^d$ Via an Affine Random
Process}
\author{Martin Hildebrand\\Department of Mathematics and Statistics\\
University at Albany\\State University of New York\\Albany, NY 12222
\\ \\ Joseph McCollum\\Department of Mathematical Sciences\\
Elms College\\Chicopee, MA 01013}

\maketitle

\begin{abstract}
This paper considers some random processes of the form
$\vcxno=T\vcxn+\vcbn \pmod p$ where $\vcbn$ and $\vcxn$ are random
variables over $(\grpz/p\grpz)^d$ and $T$ is a fixed
$d$x$d$ integer matrix which is invertible over the complex
numbers. For a particular distribution for $\vcbn$, this paper
improves results of Asci to show that if $T$ has no complex eigenvalues of
length $1$, then for integers $p$ relatively prime
to $\det(T)$, order $(\log p)^2$ steps suffice to make $\vcxn$
close to uniformly distributed where $\vcxz$ is the zero vector.
This paper also shows that if
$T$ has a complex eigenvalue which is a root of unity, then
order $p^b$ steps are needed for $\vcxn$ to get close to uniformly distributed
for some positive value $b\le 2$
which may depend on $T$ and $\vcxz$ is the zero vector.

\end{abstract}

\section{Introduction}

Previous work has looked at the following random processes on $\grpz/p\grpz$.
These processes were of the form
\[
X_{n+1}=a_nX_n+b_n\pmod p
\]
where $X_0=0$ and $a_n$ and $b_n$ had certain probability distributions.
See \cite{cdg}, \cite{mvh}, \cite{mvhtwo}, and \cite{mvhthree}. The study of 
these random processes was inspired by some pseudorandom number generators
used by computers; see, for example, \cite{knuth}. For many
choices of the distributions for $a_n$ and $b_n$, the distribution for
$X_n$ would be close to uniformly distributed for relatively small values
of $n$ (e.g. order $(\log p)^2$ or order $(\log p)\log(\log p)$)
provided that $p$ is chosen to avoid ``parity'' problems. For some choices
of $a_n$ (e.g. $a_n=1$ always), $n$ would have to be much larger
(e.g. order $p^2$) for the distribution for $X_n$ to be close to
uniformly distributed.

One can examine similar random processes on other finite structures.
For example, Asci~\cite{asci} considered random processes of the form
\[
\vcxno=T\vcxn+\vcbn \pmod p
\]
where 
$\vcxn, \vcbn$ are random variables on $(\grpz/p\grpz)^d$ with
$\vcbz, \vcbo, \dots$ i.i.d., $d$ is constant, $\vcxz=\vcz$, and $T$
is a fixed $d$x$d$ integer
matrix. The main results of \cite{asci} assumed that $T$
had non-zero integer eigenvalues. One result in \cite{asci} does not have
such an assumption but only requires $n$ to be of order $p^2(\log p)$
for $\vcxn$ to be close to uniformly distributed on $(\grpz/p\grpz)^d$.
This paper will give significantly better upper bounds provided that the
eigenvalues of $T$ are non-zero and all have length different than $1$. Such 
bounds apply, for example, to the matrix
\[
\pmatrix{2&1\cr 1&1},
\]
which does not have integer eigenvalues.
In this paper, we shall consider a specific distribution for $\vcbn$, but
the arguments potentially can be extended to other distributions for
$\vcbn$.

\section{Definitions and Main Results}

Throughout this paper, we shall assume $d$ is constant and $T$ is a fixed
$d$x$d$ integer matrix with non-zero determinant. We assume $\vcxz=\vcz$ and
$P(\vcbn=\vcz)=P(\vcbn=\vcea)=\dots=P(\vcbn=\vced)=1/(d+1)$ where 
$\vcek$ is the vector in
$(\grpz/p\grpz)^d$ whose $k$-th coordinate is $1$ and whose
other coordinates are $0$. We assume $\vcbz, \vcbo, \dots $ are i.i.d. and
\[
\vcxno=T\vcxn+\vcbn {\pmod p}.
\]
We let $P_n(s)=\Pr(\vcxn=s)$. Recall the variation distance
\begin{eqnarray*}
\|P_n-U\|&=&{1 \over 2}\sum_{s\in G}|P_n(s)-1/|G||\\
&=&\max_{A\subseteq G}|P_n(A)-U(A)|
\end{eqnarray*}
where $G$ is a finite group (here $G=(\grpz/p\grpz)^d$) and $U$ is
the uniform distribution on $G$ (i.e. $U(s)=1/|G|$ for all $s\in G$).

One of our results is the following theorem.
\begin{theorem}
\label{nolengthone}
Suppose $T$ has no eigenvalues of length $1$ over $\grpc$. 
For some value $C>0$ not
depending on $p$, 
if $n\ge C(\log p)^2$, then
$\|P_n-U\|\rightarrow 0$ as $p\rightarrow \infty$ provided
that $p$ is restricted to integers which are relatively prime to $\det(T)$.
\end{theorem}
Note that $p$ need not be prime.

Another result deals with some cases where an eigenvalue of $T$ has
length $1$.

\begin{theorem}
\label{somelengthone}
Suppose that $T$ has an eigenvalue which is a root of unity over
$\grpc$. There exists a positive value $b\le 2$ not depending on $p$ such
that if $n\le p^b$, then $\|P_n-U\|\rightarrow 1$ as $p\rightarrow \infty$
provided that $p$ is restricted to the prime numbers. 
\end{theorem}

\section{Background for Proofs of the Theorems}

The proofs of these theorems will involve the Fourier transform of
$P_n$. A more extensive description of this area appears in
Diaconis~\cite{diaconis}. A representation $\rho$ on a finite group $G$
is a map from $G$ to $GL_n(\grpc)$ such that $\rho(s)\rho(t)=
\rho(st)$ for all $s, t\in G$. The value $n$ is called the degree of the
representation and is denoted $d_{\rho}$. A representation $\rho$ with degree
$n$ is said to be irreducible if whenever $W$ is a subspace of $\grpc^n$
with $\rho(s)W\subseteq W$ for all $s\in G$, either $W=\{0\}$ or $W=\grpc^n$.
The trivial representation is given by $\rho(s)=(1)$ for all $s\in G$.
If $\rho_2(s)=M\rho_1(s)M^{-1}$ for some invertible complex matrix $M$, then
the representations $\rho_1$ and $\rho_2$ are said to be equivalent. Every 
irreducible representation on a finite group $G$ is equivalent to a unitary 
representation, i.e. a representation $\rho$ such that $(\rho(s))^{-1}=
(\rho(s))^{*}$ for all $s\in G$
where $(\rho(s))^{*}$ is the conjugate transpose of $\rho(s)$.
If $\rho$ is an irreducible representation 
on a finite group $G$ and $P$ is a probability on $G$, we define 
the Fourier transform 
\[
\hat P(\rho)=\sum_{s\in G}P(s)\rho(s).
\]

The irreducible representations of $(\grpz/p\grpz)^d$ all have degree
$1$ and are given by
\[
\rho_{\vcc}(\vcb):=q^{\sum_{i=1}^db_ic_i}
\]
where $q:=q(p):=e^{2\pi i/p}$ and
$\vcb:=\pmatrix{b_1\cr b_2\cr \vdots\cr b_d}$ and
$\vcc:=\pmatrix{c_1\cr c_2\cr \vdots\cr c_d}$ with $b_i, c_i\in \grpz/p\grpz$.
The representation $\rho_{\vcc}$ where $\vcc$ is the zero vector is the 
trivial representation.

We shall use the following lemma of Diaconis and Shahshahani. This lemma,
known as the Upper Bound Lemma, is proved in \cite{diaconis}.
\begin{lemma}
\label{upperboundlemma}
Let $P$ be a probability on a finite group $G$, and let $U$ be the uniform
distribution on $G$. Then
\[
\|P-U\|^2\le {1 \over 4}\sum_{\rho}^{*}d_{\rho}\Tr(\hat P(\rho)
\hat P(\rho)^{*})
\]
where the sum is over non-trivial irreducible representations $\rho$ such
that $\rho$ is unitary and exactly one member of each equivalence class
of non-trivial irreducible representations is included in the sum, 
$d_{\rho}$ is the degree of $\rho$, and $*$ of a matrix
denotes its conjugate transpose.
\end{lemma}

\section{Proof of Theorem~\protect\ref{nolengthone}}

First we shall develop a recurrence relation relating the Fourier
transform of $P_{n+1}$ to the Fourier transforms of $P_n$ and
$P_1$. In doing so,
we shall let $\hat P_n(\vcc)$ denote $\hat P_n(\rho_{\vcc})$. The relation
is given by the following lemma, due to Asci~\cite{asci}.

\begin{lemma}
\label{recurrencerelation}
If $p$ and $\det(T)$ are relatively prime, then
\[
\hat P_{n+1}(\vcc)={1\over d+1}\left(\hat P_n(T^t\vcc)\right)\left(
1+\sum_{r=1}^dq^{c_r}\right)
\]
where $T^t$ is the transpose of $T$.
\end{lemma}

We shall show that if 
$\vcc\ne\vcz$,
then in the sequence $\vcc$, $T^t\vcc$, $(T^t)^2\vcc$, $\dots$, 
at least one of the
coordinates of at least 
one of the first $C_2\log p$ terms (where $C_2$ is a constant)
will not be ``near'' $0$ mod $p$.

\begin{lemma}
\label{notnearzero}
Suppose $T$ has no eigenvalues of length $1$ over $\grpc$.
Suppose $\vcc\ne\vcz$. Write $\vcc\in (\grpz/p\grpz)^d$ as an element of
$\grpz^d$ with its entries as close to $0$ as possible (that is
$\vcc\in[-p/2,p/2]^d$). Then for some
positive values $C_1$ and $C_2$ (depending on $T$ but not $p$ or $\vcc$),
for sufficiently large $p$, $(T^t)^{\ell}\vcc$ has a coordinate
(viewed in $\grpz^d$) of length at least $C_1p$ but no more than $p/2$
for some non-negative $\ell\le C_2\log p$.
\end{lemma}

{\bf Proof:} Note that since $T$ is invertible over $\grpc$, then so is $T^t$.
Thus, since $\vcc\ne\vcz$, we conclude $(T^t)\vcc\ne\vcz$, 
$(T^t)^2\vcc\ne\vcz$, etc.

Let's write $T^t$ in Jordan block diagonal form over $\grpc$:
\[
T^t=M^{-1}
\pmatrix{J_1&0&\cdots\cr
0&J_2&\cdots\cr
\vdots&\vdots&\ddots}M
\]
for some invertible complex matrix $M$ where
\[
J_i=
\pmatrix{a_i&1&0&\cdots&0\cr
0&a_i&1&\cdots&0\cr
\vdots&\vdots&\vdots&\ddots&\vdots\cr
0&0&0&\cdots&a_i}.
\]
By the assumption on the eigenvalues of $T$, $|a_i|\ne 1$. Since $T$
is invertible over $\grpc$, $a_i\ne 0$.
We can write
\[
(T^t)^{\ell}=M^{-1}
\pmatrix{J_1^{\ell}&0&\cdots\cr
0&J_2^{\ell}&\cdots\cr
\vdots&\vdots&\ddots}M.
\]

The following propositions will be used in the proof of
Lemma~\ref{notnearzero}.

\begin{proposition}
\label{propone}
There exist positive constants $C_3$ and $C_4$ such that if $\vcv\in
\grpc^d$ has $v_{\arge}$ as the largest length of a coordinate, then the
coordinate of $M^{-1}\vcv$ with largest length has length at least 
$C_3v_{\arge}$ but no more than $C_4v_{\arge}$ while the coordinate of
$M\vcv$ with largest length has length at least $(1/C_4)v_{\arge}$ but no
more than $(1/C_3)v_{\arge}$.
\end{proposition}

\begin{proposition}
\label{proptwo}
If $M^{-1}\vcv\in\grpz^d$ and $M^{-1}\vcv\ne\vcz$, then at least one
of the coordinates of $\vcv$ must have length at least $C_5$ for some
constant $C_5>0$.
\end{proposition}

\begin{proposition}
\label{propthree}
For some constant $C_6>0$, if $\ell>C_6\log p$, then all coordinates
of
\[
\pmatrix{J_1^{\ell}&0&\cdots\cr 0&J_2^{\ell}&\cdots\cr \vdots&\vdots&\ddots}
M\vcc
\]
which correspond to an eigenvalue of length less than $1$ will, 
for sufficiently large $p$, all have length less than $C_5$.
\end{proposition}

\begin{proposition}
\label{propfour}
Suppose $C_1$ is such that $(C_1/C_3)(|a|+1)<1/(2C_4)$ where $a$ is the 
eigenvalue of $T$ with the largest length. Suppose all
coordinates of $M\vcc$ have length no more than $(C_1/C_3)p$.
Suppose some coordinate of
\[
\pmatrix{J_1^{\ell^{\prime}}&0&\cdots\cr 0&J_2^{\ell^{\prime}}&\cdots
\cr \vdots&\vdots&\ddots}M\vcc
\]
which corresponds to an eigenvalue of length greater than $1$ has
length greater than $C_5$ for some
$\ell^{\prime}\le (C_6\log p)+1$. Then for some constant 
$C_2>0$, there exists a value $s<C_2\log p$ such that
for sufficiently large $p$, the coordinate with the largest length of
\[
\pmatrix{J_1^s&0&\cdots\cr 0&J_2^s&\cdots\cr \vdots&\vdots&\ddots}M\vcc
\]
has length at least $(C_1/C_3)p$ but no more than
$p/(2C_4)$.
\end{proposition}

{\bf Proof of Proposition~\ref{propone}:} Let $A$ be the largest
length of the entries of $M^{-1}$. Then all coordinates of $M^{-1}\vcv$
will have
length at most $Adv_{\arge}$. Now let $B$ be the largest length of the
entries of $M$. The largest coordinate of $M^{-1}\vcv$ must have length at
least $(1/(Bd))v_{\arge}$; otherwise all coordinates of
$\vcv=M(M^{-1}\vcv)$ would have length under $v_{\arge}$. The statement
about $M\vcv$ follows directly. \qed

{\bf Proof of Proposition~\ref{proptwo}:} Let $C_5=1/C_4$, and note that 
$M^{-1}\vcv$ 
must have at least one coordinate of length at least $1$. Then use
Proposition~\ref{propone}. \qed

To prove Propositions \ref{propthree} and \ref{propfour}, we shall use the
following proposition.

\begin{proposition}
\label{blockpower}
If
\[
J=\pmatrix{a&1&0&\cdots&0\cr 0&a&1&\dots&0\cr
\vdots&\vdots&\vdots&\ddots&\vdots\cr 0&0&0&\cdots&a},
\]
then
\[
(J^{\ell})_{ij}=
\cases{a^{\ell}&if $i=j$\cr   {\ell \choose j-i}a^{\ell-(j-i)}& if $i<j$\cr
0&if $i>j$.}
\]
\end{proposition}
Note that ${m \choose n}$ is $0$ by convention if $m$ and $n$ are non-negative
integers with $m<n$.

{\bf Proof:} We shall proceed by induction on $\ell$. Note that the result is
true if $\ell=0$ or $\ell=1$. Suppose $J^{\ell}$ satisfies the
proposition. We wish to show that $J^{\ell+1}$ satisfies the analogous
result.

If $i>j$, note that since $J^{\ell}$ is an upper triangular matrix
by the induction hypothesis and $J$ is also an upper triangular
matrix, then $J^{\ell+1}$ is upper triangular and
$(J^{\ell+1})_{ij}=0$.

If $i=j$, then $(J^{\ell+1})_{ij}=\sum_k (J^{\ell})_{ik}J_{ki}=
(J^{\ell})_{ii}J_{ii}=a^{\ell}a=a^{\ell+1}$.

If $i<j$, then \begin{eqnarray*}
(J^{\ell+1})_{ij}&=&\sum_k(J^{\ell})_{ik}J_{kj}\\
&=&(J^{\ell})_{i,j-1}J_{j-1,j}+(J^{\ell})_{ij}J_{jj}\\
&=&{\ell\choose j-1-i}a^{\ell-(j-1-i)}+{\ell\choose j-i}a^{\ell-(j-i)}a\\
&=&\left({\ell \choose j-i-1}+{\ell \choose j-i}\right)a^{(\ell+1)-(j-i)}\\
&=&{\ell+1\choose j-i}a^{(\ell+1)-(j-i)}.
\end{eqnarray*}
\qed

{\bf Proof of Proposition~\ref{propthree}:} Note that the largest
length of the coordinates of $M\vcc$ is at most $C_8p$ where
$C_8=1/(2C_3)>0$ is a value depending on $M$ (and hence $T$) but not $p$.

Suppose $a$ is the eigenvalue (with length less than $1$) being considered,
$J$ is a Jordan block for the
eigenvalue $a$, and $\proj(M\vcc)$ is the projection of
$M\vcc$ onto the coordinates corresponding to this Jordan
block. Then the coordinates of $J^{\ell}\proj(M\vcc)$ have length at most
\begin{eqnarray*}
d\max\left(|a|^{\ell},{\ell\choose 1}|a|^{\ell-1},\dots,
{\ell\choose d-1}|a|^{\ell-d+1}\right)C_8p&\le&d\ell^d|a|^{\ell-d+1}C_8p\\
&\rightarrow&0
\end{eqnarray*}
as $p\rightarrow\infty$ if $(\log |a|)C_6<-1$ and $\ell>C_6\log p$. The 
proposition follows.

\qed

{\bf Proof of Proposition~\ref{propfour}:} In a Jordan block of size $1$ and
corresponding eigenvalue $b$ with $|b|>1$, the 
result should be straightforward since
$C_5|b|^m\ge (C_1/C_3)p$ for some $m$ no more than a multiple of $\log p$ and
since $(C_1/C_3)p|b|<p/(2C_4)$.

In a Jordan block $J$ of size $c$ corresponding to eigenvalue $b$ with 
$|b|>1$, let's consider the coordinates of
\[
J^{\ell^{\prime}}\proj(M\vcc)=
\pmatrix{b&1&0&\cdots&0\cr 0&b&1&\cdots&0\cr 0&0&b&\cdots&0\cr
\vdots&\vdots&\vdots&\ddots&\vdots\cr 0&0&0&\cdots&b}^{\ell^{\prime}}\proj(M\vcc)
\]
where $\proj(M\vcc)$ is the projection of $M\vcc$ onto the coordinates 
corresponding to this Jordan block. 
Let $\vcy=\pmatrix{y_1\cr \vdots \cr y_c}=J^{\ell^{\prime}}\proj(M\vcc)$.
If $|y_i|\ge C_5$ for some $i\in\{1,\dots,c\}$, then at least one of the
following statements fails:

Statement $1$: $|y_c|<C_5/p^{c-1}$.

Statement $2$: $|y_{c-1}|<C_5/p^{c-2}$.

etc.

Statement $c$: $|y_1|<C_5$.

Suppose statement $e$ is the first of these statements to fail.
Suppose $r\ge 0$. Let $\vcw=\pmatrix{w_1\cr \vdots \cr w_c}=J^{\ell^{\prime}
+r}\proj(M\vcc)=J^r\vcy$. Then 
\[
w_{c-e+1}=b^ry_{c-e+1}+{r \choose 1}b^{r-1}y_{c-e+2}+\dots+
{r\choose e-1}b^{r-(e-1)}y_c
\]
and
\[
|w_{c-e+1}|\ge {C_5\over p^{c-e}}|b|^r-\sum_{i=1}^{e-1}{r\choose i}|b|^{r-i}
{C_5\over p^{c-e+i}}
\ge{C_5|b|^r\over p^{c-1}}\left(1-{cr^c\over |b|p}\right),
\]
which for sufficiently large $p$
is at least $(C_1/C_3)p$ for some value of $r$ no more than a multiple
of $\log p$. (Choose $r=\lfloor(c+1)\log p/\log |b|\rfloor$.)
This multiple depends on $c$ and $b$ but does not depend on $p$.
However, for a given matrix $T$, there are finitely many choices for $c$
and $b$. So for some value $s$ no more than a multiple of $\log p$
(where the multiple only depends on $T$), some coordinate of
$J^s\proj(M\vcc)$ has length
at least $(C_1/C_3)p$ provided that $p$ is sufficiently 
large. Furthermore, if this value $s$ is the smallest non-negative 
value for which this statement  holds, then this coordinate has length
no more
than $(1/2C_4)p$ since $(C_1/C_3)(|b|+1)<1/(2C_4)$.
We can suppose that this coordinate is the coordinate with the
largest length.

\qed

To complete the proof of Lemma~\ref{notnearzero}, note that if some 
coordinate of $\vcc$ has length at least $C_1p$ (where $C_1$ is defined
in Proposition~\ref{propfour}), then the lemma follows
directly with $\ell=0$. Otherwise all entries in $\vcc$ have length less than $C_1p$. Thus
by Proposition~\ref{propone}, all entries of $M\vcc$ must have length less
than $(C_1/C_3)p$. By Proposition~\ref{propthree}, if $\ell>C_6\log p$, 
then all coordinates of
\[
\pmatrix{J_1^{\ell}&0&\cdots\cr 0&J_2^{\ell}&\cdots\cr \vdots&\vdots&\ddots}
M\vcc
\]
which correspond to an eigenvalue of length less than $1$ will, for
sufficiently large $p$, all have length less than $C_5$. By 
Proposition~\ref{proptwo}, some coordinate of 
\[
\pmatrix{J_1^{\ell}&0&\cdots\cr 0&J_2^{\ell}&\cdots\cr \vdots&\vdots&\ddots}
M\vcc
\]
will have length at least $C_5$. Thus by Proposition~\ref{propfour}, there
exists a value $s<C_2\log p$ such that for sufficiently large $p$,
the coordinate with the largest length of
\[
\pmatrix{J_1^s&0&\cdots\cr 0&J_2^s&\cdots\cr \vdots&\vdots&\ddots}M\vcc
\]
has length at least $(C_1/C_3)p$ but no more than $p/(2C_4)$,
and so by Proposition~\ref{propone}, 
\[
(T^t)^s\vcc=M^{-1}
\pmatrix{J_1^s&0&\cdots\cr 0&J_2^s&\cdots\cr \vdots&\vdots&\ddots}M\vcc
\]
has a coordinate of length at least $C_1p$ but no more than $p/2$ if
$p$ is sufficiently large.
\qed

To prove Theorem~\ref{nolengthone}, note that if
$C_1p\le |c_r|\le p/2$ for some $r$ in $\{1,\dots,d\}$, then
\[
\left|1+\sum_{r=1}^dq^{c_r}\right|\le C_9
\]
 for some constant
$C_9<d+1$ and
\[
|\hat P_{1}(\vcc)|={\left|1+\sum_{r=1}^dq^{c_r}\right|\over d+1}\le 
{C_9\over d+1}<1;
\]
otherwise if $\vcc\ne\vcz$, then
\[
|\hat P_{1}(\vcc)|\le 1.\]
Suppose $p$ is sufficiently large and
$s=\lfloor C_2\log p\rfloor$. Lemma~\ref{recurrencerelation},
Lemma~\ref{notnearzero}, and the
previous observation imply
\[
|\hat P_{s+1}(\vcc)|=\prod_{j=0}^s|\hat P_{1}((T^t)^j\vcc)|\le {C_9\over d+1}
\]
for all $\vcc\ne\vcz$. Thus if $r$ is a positive integer and $\vcc\ne\vcz$, then
\[
|\hat P_{r(s+1)}(\vcc)|=\prod_{i=0}^{r-1}\prod_{j=0}^s
|\hat P_{1}((T^t)^{j+i(s+1)}\vcc)|\le\left({C_9\over d+1}\right)^r.
\]
If $r$ is a sufficiently large multiple of $\log p$, then
\[
\left(C_9\over d+1\right)^r<{1\over p^{d+1}},
\]
and 
\[
\sum_{\vcc\ne\vcz}|\hat P_{r(s+1)}(\vcc)|^2<{1\over p^2}.
\]
Note that $r(s+1)<C(\log p)^2$ for some constant $C>0$.
Then, from Lemma~\ref{upperboundlemma} and Lemma~\ref{recurrencerelation},
if $n\ge C(\log p)^2$, then
\[
\|P_n-U\|\le \|P_{r(s+1)}-U\|\le
{1\over 2}\sqrt{\sum_{\vcc\ne\vcz}|\hat P_{r(s+1)}(\vcc)|^2}<{1\over 2p}
\rightarrow 0\]
as $p\rightarrow \infty$.
\qed

\section{Proof of Theorem~\protect\ref{somelengthone}}

Suppose $T$ has an eigenvalue (over $\grpc$) which is an $m$-th root
of unity. Then $T^m$ has an eigenvalue $1$ over
$\grpc$.
So the characteristic polynomial of $T^m$ over $\grpc$ has $1$ as a 
root. Thus the characteristic polynomial of $T^m$ over $\grpz/p\grpz$
has $1$ as a root, and $1$ is an eigenvalue of $T^m$ over
$\grpz/p\grpz$. 
We can write $(T^m)^t=MJM^{-1}$ for some matrix $J$ in Jordan block
form over $\grpz/p\grpz$ and so $T^m=(M^{-1})^tJ^tM^t$. Thus there
exists a basis $\vcba, \dots, \vcbd$ for $(\grpz/p\grpz)^d$ over
$\grpz/p\grpz$ such that if $p(\vcv)=a_1$ when 
$\vcv=a_1\vcba+\dots+a_d\vcbd$, we
have $p(T^m\vcba)=1$ and $p(T^m\vcbj)=0$ if $j=2,\dots,d$.
Note that for all $n$,
\[
\vcxnm=T^m\vcxn+T^{m-1}\vcbn+\dots+T\vcbnmt+\vcbnmo \pmod p,
\]
and so $p(\vcxnm)=p(T^m\vcxn)+p(T^{m-1}\vcbn)+\dots+p(T\vcbnmt)+p(\vcbnmo)$.
Note that $p(T^m\vcxn)=p(\vcxn)$. There are at most $d+1$ possible values
of $p(T^{m-1}\vcbn)$, $d+1$ possible values of $p(T\vcbnmt)$, and $d+1$
possible values of $p(\vcbnmo)$. So there are at most $(d+1)^m$ possible
values of $p(T^{m-1}\vcbn)+\dots+p(T\vcbnmt)+p(\vcbnmo)$.
Thus $p(\vcxz), p(\vcxm), p(\vcxtm), \dots$ forms a random walk on
$\grpz/p\grpz$ with support of size $u$ which is at most $(d+1)^m$.
By
Greenhalgh~\cite{greenhalgh}, given $\epsilon>0$, there exists a value $C_2>0$ such that if
$c<C_2p^{2/(u-1)}$, then the variation distance of $p(\vcxcm)$ from
uniform in $\grpz/p\grpz$ is at least $1-\epsilon$. Thus, in the
random process on $(\grpz/p\grpz)^d$, we get
$\|P_{cm}-U\|>1-\epsilon$. Results of Asci~\cite{asci} imply that
$b$ is no larger than $2$. The theorem follows.
\qed

\section{Questions for Further Study}

One question to consider is whether the order $(\log p)^2$
rate of convergence can be improved to order $(\log p) \log (\log p)$
similar to Theorem 2 in \cite{mvhtwo}. For diagonal matrices with non-zero
eigenvalues all different from $\pm 1$, Theorem 4.5 of Asci~\cite{asci}
confirms this, but for more general
cases, the question is open. Another question is analogous to one
explored in Chung, Diaconis, and Graham~\cite{cdg}: Are order $\log p$
steps sufficient for typical values of $p$?

Another question for further study is to see how large the value $b$ in
Theorem~\ref{somelengthone} can be. Also, given $\epsilon>0$, will there exist
a $C>0$ such that if $n\le Cp^2$, then $\|P_n-U\|>1-\epsilon$ for
sufficiently large $p$?

Another area to explore involves random processes on fields with $p^d$
elements where $p$ is a prime; these random processes would be of the
form $X_{n+1}=AX_n+b_n$ where $b_n$ and $X_n$ are random variables over this
field and $A$ is a fixed element of this field. Since elements
of this field can be represented as a $d$x$d$ matrix over $\grpz/p\grpz$
(see, for example, pp. 64-65 of \cite{ln}), perhaps some natural connections 
can be made between the random processes studied here and such random 
processes on finite fields.

\section{Acknowledgments}

The authors would like to thank Persi Diaconis for posing a related
question to one of the authors. The authors would like to thank the referee
for some comments incorporated into this paper.

\end{document}